\documentclass[12pt]{article}
\usepackage[utf8]{inputenc}
\usepackage{amsmath}
\usepackage{amsfonts}
\usepackage{amssymb}
\usepackage{amsthm}
\usepackage{relsize}
\usepackage{array}
\usepackage{multirow}
\usepackage{tabu}
\usepackage{soul}
\usepackage{graphicx}
\usepackage{epigraph}
\usepackage{float}
\usepackage[english]{babel}
\usepackage[usenames, dvipsnames]{color}
\usepackage{fancyhdr}
\usepackage[Sonny]{fncychap}
\usepackage{tikz-cd}
\usepackage{tkz-euclide}
\usepackage{sseq}
\usepackage{newfloat}
\usepackage{amsmath,mathtools}
\usepackage{pgfplots}
\usepackage{geometry}
 \usepackage{sectsty}
\sectionfont{\fontsize{14}{15}\selectfont}

\usepackage[all]{xy}
\DeclareFloatingEnvironment[fileext=dia]{diagram}

\theoremstyle{definition}

\usepackage[most]{tcolorbox}

\makeatletter
\newtcolorbox{note}[1][]{%
	breakable,
	enhanced jigsaw, 
	borderline west={3pt}{0pt}{black!10!white}, 
	borderline south={1pt}{0pt}{black!10!white}, 
	borderline east={1pt}{0pt}{black!10!white},
	borderline north={1pt}{0pt}{black!10!white},
	sharp corners, 
	boxrule=0pt, 
	attach title to upper, 
	left=0pt,
	right=0pt,
	top=0pt,
	bottom=0pt,
	boxsep=5pt,
	colback=white,
	frame hidden,
	#1
}
\makeatother

\makeatletter
\newtcolorbox{note1}[1][]{%
	breakable,
	enhanced jigsaw, 
	sharp corners, 
	boxrule=0pt, 
	attach title to upper, 
	fontupper=\linespread{1.1}\fontfamily{qpl}\selectfont,
	fontlower=\linespread{1.1}\fontfamily{qpl}\selectfont, 
	left=0pt,
	right=0pt,
	top=0pt,
	bottom=0pt,
	boxsep=3pt,
	colback=green!3!white,
	frame hidden,
	before skip=10pt plus 2pt,after skip=10pt plus 2pt,
	#1
}
\makeatother

\makeatletter
\newcommand\tabfill[1]{%
	\dimen@\linewidth
	\advance\dimen@\@totalleftmargin
	\advance\dimen@-\dimen\@curtab
	\parbox[t]\dimen@{#1\ifhmode\strut\fi}%
}
\makeatother

\usepackage{tabularx}
\usepackage{imakeidx}
\usepackage{hyperref}
\hypersetup{colorlinks={true},linkcolor={blue},citecolor={blue},urlcolor={blue}, linktoc=all}  
 
 \usepackage[capitalize, nameinlink]{cleveref}
 \crefdefaultlabelformat{#2\textbf{#1}#3}
 \crefname{figure}{Figure}{Figures} 
 \Crefname{figure}{Figure}{Figures}
 \crefname{table}{Table}{Tables}
 \Crefname{table}{Table}{Tables}
 \crefname{section}{\S\hspace{-1mm}}{\S\hspace{-1mm}}
 \Crefname{section}{\S\hspace{-1mm}}{\S\hspace{-1mm}}
 \crefname{equation}{}{}
 \Crefname{equation}{}{}
 \crefname{example}{Geometric Pattern}{Geometric Patterns} 
 \Crefname{example}{Geometric Pattern}{Geometric Patterns}

 \usepackage{graphicx,tipa}

 \usepackage{footnote}

\begin{document}

\title{\textbf{Excavation Problems in Elamite Mathematics}}

\author{Nasser Heydari\footnote{Email: nasser.heydari@mun.ca}~ and  Kazuo Muroi\footnote{Email: edubakazuo@ac.auone-net.jp}}

\maketitle

\begin{abstract}
In this article, we study the  problems found in    the Susa Mathematical Texts No.\,24 and No.\,25 (\textbf{SMT No.\,24} and \textbf{SMT No.\,25}) which concern  excavation projects such as  canals and holes. We also  examine certain      Elamite structures, such as the canal systems serving Susa and a reservoir at the ziggurat of Chogha Zanbil, in whose construction  geometry might well have   played an important role.  
\end{abstract}

\section{Introduction}
\textbf{SMT No.\,24} and \textbf{SMT No.\,25} are two of  the texts inscribed on  26  clay tablets excavated from Susa in  southwest Iran by French archaeologists in 1933. The texts of all the Susa mathematical texts (\textbf{SMT}) along with their interpretations were first published in 1961 (see \cite{BR61}).

  \textbf{SMT No.\,24} and  \textbf{SMT No.\,25} are on a   one-column clay tablet\footnote{The reader can see the new  photos of this tablet on the website of the Louvre's collection. Please see \url{https://collections.louvre.fr/en/ark:/53355/cl010186434} for obverse  and   reverse.}. Following Bruins in \cite{BR61}, we treat these two texts  separately.

  \textbf{SMT No.\,24}  contains two problems. The first, an excavation problem, which leads to a complicated quadratic equation, is found on the obverse of the tablet. The second problem,  concerning the digging of a canal,  is presented on the reverse.    

The text of  \textbf{SMT No.\,25},  also on the reverse of the tablet, is another belonging to the category of     excavation problems. There is only one problem  treated  in this text, which concerns    digging  a canal.  
 
 Although the entire problems are unavailable because of damage to the tablet, considerable understanding of the mathematical calculations utilized in solving these problems can be derived from a careful analysis of the text that remains.
 
\section{Excavation Problems in the SMT}
In Babylonian mathematics, there are many problems which concern digging a hole, a cistern, a canal or a foundation for a building or a wall. Such problems are often referred as \textit{excavation problems}. Among the  Babylonian mathematical texts,  there are several    which address  various quantities concerning canals, although from the mathematical point of view they are relatively simple. By way of examples of such texts, we   mention \textbf{YBC 4666}, \textbf{YBC 7164}, \textbf{YBC 9874}, \textbf{VAT 7528}, and \textbf{BM 85196} (please also see \cite{Mur92-4}, for more details about the  texts of such tablets).


 In the following two parts, we discuss   \textbf{SMT No.\,24} and  \textbf{SMT No.\,25} and provide our mathematical interpretation of each text. 
  
\subsection{\textbf{SMT No.\,24}}

\subsubsection*{Transliteration}\label{SS-TI-SMT24}  

\begin{note1} 
	\underline{Obverse,  Lines 1-40}\\
	(L1)\hspace{2mm} [pa$_5$ \textit{i-na} g]i\v{s}-gi \textit{\`{u}} gi\v{s}-tir \textit{ak-\c{s}\'{u}-ru} [$\cdots$ $\cdots$]
	\\
	(L2)\hspace{2mm} [$\cdots$ $\cdots$] 30 dirig dagal ki-ta \textit{\v{s}\`{a}-lu-u\v{s}-ti} [$\cdots$ $\cdots$]
	\\
	(L3)\hspace{2mm} [\textit{ak-\c{s}}]\textit{\'{u}-ru \v{s}u-up-lu} 24 sahar \textit{ak-\c{s}\'{u}-r}[\textit{u} $\cdots$]
	\\
	(L4)\hspace{2mm} [\textit{\`{u}}] \textit{\v{s}u-up-lu mi-nu} za-e 1 dagal \textit{\v{s}\`{a} la t}[\textit{i-du-\'{u}}  $\cdots$]
	\\
	(L5)\hspace{2mm} 30 \textit{ta-mar} $\cdots$    [$\cdots$] 30 dirig gar 1 a-r\'{a} [$\cdots$ $\cdots$]
	\\
	(L6)\hspace{2mm} 30 dirig \textit{\v{s}\`{a}-lu-u\v{s}-ti} 30 \textit{le-q\'{e}} 10 \textit{ta-mar} [$\cdots$ $\cdots$]
	\\
	(L7)\hspace{2mm} \textit{i-\v{s}\'{i}-ma} 2 \textit{ta-mar} 2 \textit{ki-ma} sag gar 1 dagal a[n-ta   $\cdots$]
	\\
	(L8)\hspace{2mm} 1,30 \textit{ta-mar} 1/2 1,30 [\textit{he-pe}] 45 \textit{ta-mar} 45 [u\v{s} $\cdots$ $\cdots$]
	\\
	(L9)\hspace{2mm} \textit{tu-\'{u}r-ma \v{s}\`{a}-lu-}[\textit{u\v{s}-ti} 30 di]rig \textit{le-q\'{e}} 10 \textit{ta-mar} [$\cdots$ $\cdots$]
	\\
	(L10)\hspace{0mm} \textit{i-\v{s}\'{i}-ma} 2  \textit{ta-mar}   $\cdots$ [$\cdots$ $\cdots$]
	\\
	(L11)\hspace{0mm} 15 \textit{ta-mar} 15   $\cdots$ [$\cdots$   $\cdots$ $\cdots$] u\v{s} gar   $\cdots$ [$\cdots$]
	\\
	(L12)\hspace{0mm} \textit{re-i\v{s}-ka li-ki-}[\textit{il} $\cdots$ $\cdots$]   $\cdots$ igi-45 u\v{s} \textit{pu-\c{t}\'{u}-}[\textit{\'{u}r}]
	\\
	(L13)\hspace{0mm} \textit{a-na} [$\cdots$   $\cdots$ $\cdots$]   $\cdots$ a-\v{s}\`{a} \textit{i-\v{s}\'{i}-ma} 9,22,30 [$\cdots$] \textit{i-\v{s}\'{i}}
	\\
	(L14)\hspace{0mm} [$\cdots$ $\cdots$ $\cdots$ $\cdots$   $\cdots$ 2]4(?) \textit{\`{u}} 2 \textit{ma-na-at}
	\\
	(L15)\hspace{0mm} [$\cdots$ $\cdots$ $\cdots$ $\cdots$   $\cdots$] sag(?) 2 zi a-\v{s}\`{a}
	\\
	(L16)\hspace{0mm} [$\cdots$ $\cdots$ $\cdots$ $\cdots$ \textit{i-\v{s}}]\textit{\'{i}-ma} u\v{s} 15 \textit{wa-\c{s}\'{i}-ib} u\v{s}
	\\
	(L17)\hspace{0mm} [$\cdots$ $\cdots$ $\cdots$ $\cdots$ 3]0 [\textit{a}]\textit{-na} 9,22,30 \textit{i-\v{s}\'{i}-ma} 
	\\
	(L18)\hspace{0mm} [4,41,15 \textit{ta-ma}]\textit{r} \textit{tu-\'{u}r-ma} 45 u\v{s} \textit{a-na} 2 sag \\
	(L19)\hspace{0mm} [\textit{i-\v{s}\'{i}-ma} 1,30] \textit{a-na} [9,2]2,30 \textit{i-\v{s}\'{i}-ma} 14:3,45 \textit{ta-mar}
	\\
	(L20)\hspace{0mm} [14:3,45 a-\v{s}\`{a} \textit{s\`{a}-a}]\textit{r-ru} [14:]3,45 [\textit{a-n}]\textit{a} 4,[4]1,15 a-\v{s}\`{a}
	\\
	(L21)\hspace{0mm} [\textit{i-\v{s}\'{i}-ma} 1,5,55],4,41,15 \textit{ta-mar re-}[\textit{i\v{s}-k}]\textit{a li-ki-il}
	\\
	(L22)\hspace{0mm} [15 \textit{wa-\c{s}\'{i}-ib} u\v{s} \textit{a-na} 2] sag \textit{i-\v{s}\'{i}-ma} 30 \textit{t}[\textit{a-mar} 30 \textit{n}]\textit{a-s\'{i}-ih} sag
	\\
	(L23)\hspace{0mm} [30 \textit{a-na} 3 \textit{i-\v{s}}]\textit{í-ma} [1],30 \textit{ta-mar} [3]0 \textit{i-na} 1,[30] zi
	\\
	(L24)\hspace{0mm} 1 \textit{ta-mar} 1 \textit{a-}[\textit{na}] 9,22,30 \textit{i-\v{s}\'{i}-ma} 9,2[2,30] \textit{ta-mar}
	\\
	(L25)\hspace{0mm} \textit{a\v{s}-\v{s}um} $ < $1$ > $ \textit{ki-ma} u\v{s} \textit{qa-bu-ku} 1 a-r\'{a} [\textit{a-na} 9],22,[30 dah]
	\\
	(L26)\hspace{0mm} 1,9,22,30 \textit{ta-mar} 1/2 1,[9,2]2,30 \textit{he-pe} 34,41,15 \textit{ta-mar}
	\\
	(L27)\hspace{0mm} 34,41,15 nigin 20:3,[13,21,3]3,45 \textit{ta-mar}
	\\
	(L28)\hspace{0mm} \textit{a-na} 20:3,13,21,33,[45]  1,5,55,4,41,15 dah
	\\
	(L29)\hspace{0mm} 21:9,8,26,15 [\textit{ta-mar}] \textit{m}[\textit{i-n}]\textit{a} \'{i}b-si 35,37,30 \'{i}b-si
	\\
	(L30)\hspace{0mm} 34,41,15 \textit{ta-ki-i}[\textit{l-ta}]-\textit{ka} [\textit{a-na} 3]5,37,30 dah
	\\
	(L31)\hspace{0mm} 1,10:18,45 \textit{ta-m}[\textit{ar}] \textit{m}[\textit{i-na a-na}] 14:3,45
	\\
	(L32)\hspace{0mm} a-\v{s}\`{a} \textit{s\`{a}-ar-ri} gar \textit{\v{s}\`{a}} [1,10:1]8,45 \textit{i-na-ad-di-na}
	\\
	(L33)\hspace{0mm} 5 gar 5 dagal an-ta 1/2 [5 \textit{he-pe} 2,30 \textit{ta-mar} 2,30 \textit{a-na}]
	\\
	(L34)\hspace{0mm} 30 dirig dah 3 \textit{ta-mar} 3 dagal ki-ta [igi-12 \textit{\v{s}\`{a}} dagal an-ta]
	\\
	(L35)\hspace{0mm} ugu dagal ki-ta \textit{i-te-ru} \textit{le-q}[\textit{\'{e}}] 10 \textit{ta-mar} [30 \textit{\`{u}} 10 ul-gar]
	\\
	(L36)\hspace{0mm} 40 \textit{a-na} 12 \textit{\v{s}u-up-li i-\v{s}\'{i}} 8 \textit{ta-mar} [5 dagal an-ta]
	\\
	(L37)\hspace{0mm} \textit{\`{u}} 3 dagal ki-ta ul-gar 8 \textit{ta-mar} [1/2 8 \textit{he-pe}]
	\\
	(L38)\hspace{0mm} 4 \textit{ta-mar} 4 \textit{a-na} 8 \textit{\v{s}u-up-li i-\v{s}\'{i}-}[\textit{ma}]
	\\
	(L39)\hspace{0mm} 32 \textit{ta-mar} igi-32 \textit{pu-\c{t}\'{u}-\'{u}r} 1,5[2,30 \textit{ta-mar}]\\
	(L40)\hspace{0mm} 1,52,30 \textit{a-na} 24 sahar \textit{i-\v{s}}[\textit{\'{i}} 45 \textit{ta-mar} 45 u\v{s}]\\
	
	\underline{Reverse,  Lines 1-24}\\
	(L1)\hspace{2mm} za-e(?) $\cdots$ [$\cdots$ $\cdots$ $\cdots$]
	\\
	(L2)\hspace{2mm} \textit{\'{u}-\v{s}\`{a}-p\'{i}-i}[\textit{l}] \textit{i}(?)-[\textit{na}(?) \textit{k}]\textit{a-la-ak-}[\textit{ki-im} $\cdots$ $\cdots$]
	\\
	(L3)\hspace{2mm} 2-kam \textit{ta-a}[\textit{d}]-\textit{di-in} 2-\textit{tu} [$\cdots$ $\cdots$]
	\\
	(L4)\hspace{2mm} a-\v{s}\`{a} \textit{ka-la-ak-ki} gal ul-[gar $\cdots$ $\cdots$]
	\\
	(L5)\hspace{2mm} \textit{a-na} t\`{u}n \textit{\v{s}\`{a} a-ta-ap pa}-[\textit{n}]\textit{a-n}[\textit{im}] da[h $\cdots$]
	\\
	(L6)\hspace{2mm} \textit{a-na} tur \textit{a\v{s}}-[\textit{lu-u}]\textit{\c{t}} (?) ul-gar sahar $\cdots$ [$\cdots$ $\cdots$ 1,15]
	\\
	(L7)\hspace{2mm} za-e 1,15 ul-gar \textit{a-na} 13 \textit{\v{s}\`{a}}-[\textit{la-\v{s}\'{e}-ra-t}]\textit{i} \textit{i-\v{s}\'{i}-ma} 16,[15]
	\\
	(L8)\hspace{2mm} 10 \textit{\v{s}\`{a} ka-la-ak-ku} ugu \textit{ka-la-a}[\textit{k-ki i-t}]\textit{e}-\textit{ru} nigin 1,40 \textit{ta-mar}
	\\
	(L9)\hspace{2mm} 1,40 \textit{i-na} 16,15 zi 16,[13,20] \textit{ta-mar} igi-10 dirig \textit{pu-ţ\'{u}-\'{u}r}
	\\
	(L10)\hspace{0mm} 6 \textit{ta-mar} igi-12 \textit{\v{s}u-up-li pu-ţ\'{u}-\'{u}r} 5 \textit{ta-mar} 5 \textit{a-na} 6 \textit{i-\v{s}\'{i}}
	\\
	(L11)\hspace{0mm} 30 \textit{ta-mar}  30 \textit{ta-lu-ku} 30 \textit{ta-lu-ka a-na} 16,13,20 \textit{i-\v{s}\'{i}-ma}
	\\
	(L12)\hspace{0mm} 8,6,40 \textit{ta-mar} 10 [dir]ig nigin 1,40 \textit{ta-mar} 1,40 \textit{a-na} 13 \textit{\v{s}\`{a}-la-\v{s}\'{e}-ra-ti}
	\\
	(L13)\hspace{0mm} \textit{i-\v{s}\'{i}-ma} 21,[40 \textit{ta-mar}] 21,40 \textit{i-na} 8,6,40 zi
	\\
	(L14)\hspace{0mm} 7,45 \textit{ta}-[\textit{mar} \textit{re-i\v{s}-k}]\textit{a} \textit{li-ki-il} 3[0 \textit{ta-lu-k}]\textit{a}
	\\
	(L15)\hspace{0mm} \textit{a-na} 13 [\textit{\v{s}\`{a}-la-\v{s}\'{e}-ra-ti}] \textit{i-\v{s}\'{i}} 6,30 \textit{t}[\textit{a-mar} 30 \textit{t}]\textit{a-lu-ka}
	\\
	(L16)\hspace{0mm} \textit{a-na} \textit{ka-aiia-ma}-[\textit{ni}] 2 tab-ba 1 \textit{ta-mar} 1 \textit{a-na} 6,30 dah
	\\
	(L17)\hspace{0mm} [7],30 \textit{ta-mar} 1[3 \textit{\v{s}\`{a}-l}]\textit{a-a\v{s}-\v{s}\'{e}-ra-ti a-na} 3-\textit{\v{s}u a-na ka-aiia-ma-ni}
	\\
	(L18)\hspace{0mm} \textit{a-li-ik-ma} 3[9] \textit{ta-mar} 7,30 \textit{a-na} 39 dah 46,30 \textit{ta-mar}
	\\
	(L19)\hspace{0mm} \textit{mi-na a-na} 46,30 gar \textit{\v{s}\`{a}} 7,45 \textit{\v{s}\`{a}} \textit{r}[\textit{e-i\v{s}}]-\textit{ka} \textit{\'{u}-ki-il-lu}
	\\
	(L20)\hspace{0mm} [\textit{i-n}]\textit{a-ad-di-na} [10] gar \textit{re-i\v{s}-ka li-ki-il} 1/2 10 dirig \textit{he-pe} 
	\\
	(L21)\hspace{0mm} [5 \textit{ta}]-\textit{mar} [5(?) gar(?)] 5 nigin 25 \textit{ta-mar} 25 \textit{a-na} 10
	\\
	(L22)\hspace{0mm} [\textit{\v{s}\`{a} re-i\v{s}-ka \'{u}-ki-il-lu}] dah 10,25 \textit{ta-mar mi-na} \'{i}b-si 
	\\
	(L23)\hspace{0mm} [25 \'{i}b-si 25(?) gar(?)] 5 \textit{a-na} 25 \textit{i\v{s}-te-en} dah 30 \textit{ta-mar} \\
	(L24)\hspace{0mm} [\textit{i-na} 25 2-kam zi 20 \textit{t}]\textit{a-mar} 30 gal 20 tur
	
\end{note1}

\subsubsection*{Translation}\label{SS-TR-SMT24}  

\underline{Obverse,  Lines 1-40}
\begin{tabbing}
	\hspace{17mm} \= \kill 
	(L1)\> \tabfill{The canal that I constructed in reed beds and woods, $\cdots$ $\cdots$.}\index{reed bed}\index{canal}\\
	(L2)\> \tabfill{$\cdots$ 0;30 of the excess, the lower breadth. One third of $\cdots$,}\index{breadth (of a canal)}\\
	(L3)\> \tabfill{that I constructed, the depth. The volume 24,0({\fontfamily{qpl}\selectfont volume-sar}\index{volume-sar (capacity unit)}) that I constructed $\cdots$.}\index{depth (of a canal)}\\
	(L4)\> \tabfill{What are $\cdots$ and the depth? You, 1 breadth that you do not know $\cdots$,}\index{breadth (of a canal)}\index{depth (of a canal)}\\
	(L5)\> \tabfill{you see 0;30. Put down 0;30 of the excess. 1 regular number $\cdots$,}\index{normal number}\index{regular number}\\
	(L6)\> \tabfill{0;30 of the excess. Take one third of 0;30, (and) you see 0;10. $\cdots$.}\\
	(L7)\> \tabfill{Multiply (it by 12), and you see 2. Put down 2 as the width. 1, the upper breadth $\cdots$,}\index{breadth (of a canal)}\index{width}\\
	(L8)\> \tabfill{you see 1,30. Halve 1,30, (and) you see 45. 45 is the length. $\cdots$.}\index{length}\\
	(L9)\> \tabfill{Return. Take one third of 0;30 of the excess, (and) you see 0;10. $\cdots$.}\\
	(L10)\> \tabfill{Multiply (it by 12), and you see 2. $\cdots$,}\\
	(L11)\> \tabfill{you see 15. 15, $\cdots$. Put down the length. $\cdots$.}\index{length}\\
	(L12)\> \tabfill{Let it hold your head. $\cdots$. Make the reciprocal of 45 of the length, (and 0;1,20).}\index{reciprocal of a number}\index{length}\\
	(L13)\> \tabfill{Multiply (it) by $\cdots$ of the area, and you see 9;22,30. Multiply $\cdots$.}\index{area}\\
	(L14)\> \tabfill{$\cdots$ $\cdots$ 24,0(?), and 0;2 is the ratio (of the width to the length).}\index{length}\index{similarity ratio}\index{width}\\
	(L15)\> \tabfill{$\cdots$ $\cdots$ width(?), subtract 0;2(?), and the area.}\index{width}\index{area}\\
	(L16)\> \tabfill{$\cdots$ $\cdots$, and the length. 15, the one which is to be added to the length.}\index{length}\\
	(L17)\> \tabfill{$\cdots$ $\cdots$. Multiply 0;30 by 9;22,30, and}\\
	(L18-19)\> \tabfill{you see 4;41,15. Return. Multiply 45 of the length by 0;2 of (the ratio of) the width (to the length), and you see 1;30. Multiply (it) by 9;22,30, and you see 14;3,45.}\index{length}\index{similarity ratio}\index{width}\\
	(L20-21)\> \tabfill{14;3,45 is the false area. Multiply 14;3,45 by 4;41,15 of the area, and you see 1,5;55,4,41,15. Let it hold your head.}\index{provisional area}\index{area}\\ 
	(L22)\> \tabfill{Multiply 15, the one which is to be added to the length, by 0;2 of (the ratio of) the width (to the length), and you see 0;30.}\index{length}\index{similarity ratio}\index{width}\\ 
	(L23)\> \tabfill{Multiply 0;30 by 3, and you see 1;30. Subtract 0;30 from 1;30, (and)}\\ 
	(L24)\> \tabfill{you see 1. Multiply 1 by 9;22,30, and you see 9;22,30.}\\ 
	(L25)\> \tabfill{Since $ <1,0>  $ as the length is said to you, add 1,0, the factor, to 9;22,30, (and)}\index{length}\index{regular factor}\\
	(L26)\> \tabfill{you see 1,9;22,30. Halve 1,9;22,30, (and) you see 34;41,15.}\\ 
	(L27)\> \tabfill{Square 34;41,15, (and) you see 20,3;13,21,33,45.}\\ 
	(L28)\> \tabfill{Add 1,5;55,4,41,15 to 20,3;13,21,33,45, (and)}\\
	(L29)\> \tabfill{you see 21,9;8,26,15. What is the square root? 35;37,30 is the square root.}\index{square root}\\ 
	(L30)\> \tabfill{Add 34;41,15 (which was used in) your completing the square to 35;37,30, (and)}\index{completing the square method}\\
	(L31-32)\> \tabfill{you see 1,10;18,45. What should I put down to 14;3,45, the false area, which will give me 1,10;18,45?}\index{area}\\
	(L33)\> \tabfill{Put down 5. 5 is the upper breadth. Halve 5, (and) you see 2;30.}\index{breadth (of a canal)}\\ 
	(L34)\> \tabfill{Add 2;30 to 0;30 of the excess, (and) you see 3. 3 is the lower breadth.}\index{breadth (of a canal)}\\
	(L35)\> \tabfill{Take 1/12 of the amount by which the upper breadth exceeded the lower breadth, (and) you see 0;10. Add 0;30 and 0;10 together, (and the result is 0;40).}\index{breadth (of a canal)}\\
	(L36)\> \tabfill{Multiply 0;40 by 12 of the (constant of the) depth, (and) you see 8.}\index{depth (of a canal)}\\
	(L37)\> \tabfill{Add together 5 of the upper breadth and 3 of the lower breadth, (and) you see 8. Halve 8, (and)}\index{breadth (of a canal)}\\
	(L38)\> \tabfill{you see 4. Multiply 4 by 8, of the depth, and}\index{depth (of a canal)}\\
	(L39)\> \tabfill{you see 32. Make the reciprocal of 32, (and) you see 0;1,52,30.}\index{reciprocal of a number}\\
	(L40)\> \tabfill{Multiply 0;1,52,30 by 24,0 of the volume, (and) you see 45. 45 is the length.}\index{length} 		 
\end{tabbing}\index{canal}   
\noindent
\underline{Reverse,  Lines 1-24}
\begin{tabbing}
	\hspace{17mm} \= \kill 
	(L1)\> \tabfill{You(?), $\cdots$ $\cdots$ $\cdots$ $\cdots$}\\ 
	(L2)\> \tabfill{$\cdots$ $\cdots$ I excavated . In(?) the hole $\cdots$ $\cdots$.}\\ 
	(L3)\> \tabfill{you gave the second $\cdots$ $\cdots$. A second time $\cdots$ $\cdots$.}\\ 
	(L4)\> \tabfill{I added  $\cdots$ $\cdots$ and the area of the large hole together, $\cdots$ $\cdots$.}\index{area}\\ 
	(L5)\> \tabfill{I added $\cdots$ $\cdots$ to the depth\index{depth (of a canal)} of the former canal, $\cdots$ $\cdots$.}\\ 
	(L6)\> \tabfill{I cut off (?) $\cdots$ $\cdots$ for the small. The sum of the volume and $\cdots$ $\cdots$ is 1;15.}\\  
	(L7)\> \tabfill{You, multiply 1;15 of the sum by 13 of one thirteenth, and (you see) 16;15.}\\  
	(L8)\> \tabfill{Square 0;10 of the amount by which (the length of) the (large) hole exceeded (the length of) the (small) hole, (and) you see 0;1,40.}\index{length}\\  
	(L9)\> \tabfill{Subtract 0;1,40 from 16;15, (and) you see 16;13,20. Make the reciprocal of 0;10 of the excess,}\index{reciprocal of a number} \\ 
	(L10)\> \tabfill{(and) you see 6. Make the reciprocal of 12 of the depth, (and) you see 0;5. Multiply 0;5 by 6,}\index{reciprocal of a number}\index{depth (of a canal)} \\  
	(L11)\> \tabfill{(and) you see 0;30. 0;30 is the product. Multiply 0;30 of the product by 16;13,20, and}\\  
	(L12)\> \tabfill{you see 8;6,40. Square 0;10 of the excess, (and) you see 0:1,40. Multiply 0;1,40 by 13 of one thirteenth,}\\ 
	(L13)\> \tabfill{and you see 0;21,40. Subtract 0;21,40 from 8;6,40, (and)}\\  
	(L14)\> \tabfill{you see 7;45. Let it hold your head.}\\  
	(L15)\> \tabfill{Multiply 0;30 of the product by 13 of one thirteenth, (and) you see 6;30.}\\  
	(L16)\> \tabfill{Multiply 0;30 of the product by regular (number) 2, (and) you see 1. Add 1 to 6;30, (and)}\index{normal number}\index{regular number}\\  
	(L17)\> \tabfill{you see 7;30. Multiply 13 of the one thirteenth by 3, by regular (number three),}\index{normal number}\index{regular number}\\ 
	(L18)\> \tabfill{and you see 39. Add 7;30 to 39, (and) you see 46;30.}\\  
	(L19)\> \tabfill{What should I put to 46;30 which gives me 7;45 that held your head?}\\ 
	(L20)\> \tabfill{Put down 0;10. Let it hold your head. Halve 0;10 of the excess, (and)}\\ 
	(L21)\> \tabfill{you see 0;5. Put down 0;5(?). Square 0;5, (and) you see 0;0,25.}\\  
	(L22)\> \tabfill{Add 0;0,25 to 0;10 that held your head, (and) you see 0;10,25. What is the square root?}\index{square root}\\ 
	(L23)\> \tabfill{0;25 is the square root. Put down 0;25(?). On the one hand add 0;5 to 0;25, (and) you see 0;30.}\index{square root}\\  
	(L24)\> \tabfill{On the other hand subtract (0;5) from 0;25, (and) you see 0;20. 0;30 is the large, (and) 0;20 is the small.}
\end{tabbing}\index{normal number}\index{regular number}\index{canal}

\subsubsection*{Mathematical Calculations} 
\begin{figure}[H]
	\centering
	\includegraphics[scale=1]{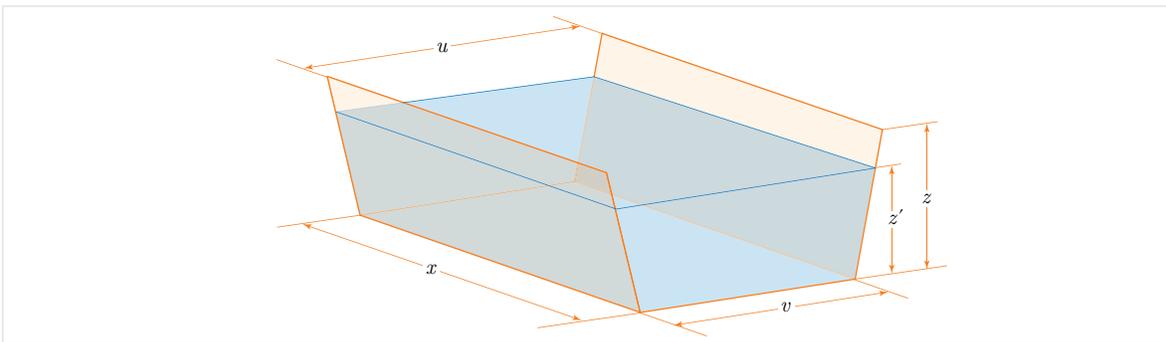}
	\caption{The general shape of a canal and its dimensions}
	\label{Figure1}
\end{figure}

The two problems in this text deal with constructing   canals  and computing their dimensions.   The general shape of  a canal  is a prism with   trapezoidal  bases as   shown in \cref{Figure1} along with its reserved water. Denote   the lower breadth (width), the upper  breadth  (width), the length and the depth  (height)  of the canal by  $v$, $u$, $x$ and $z$ respectively. Also denote the height of the reserved water by $z'$.

The area $S$ of the trapezoidal base and the volume  $V$ of the trapezoidal canal  are obtained by
\begin{equation}\label{equ-SMT24-a}
	S=\frac{1}{2}z(u+v),
\end{equation}
and
\begin{equation}\label{equ-SMT24-b}
	V=xS=\frac{1}{2}xz(u+v).
\end{equation}
Although formula \cref{equ-SMT24-b} gives us the whole capacity of a canal, it is possible to compute the capacity of its reserved water $V'$   by using the constant of a canal given in  \textbf{SMT No.\,3}, line 33. In this line we read:\\
  {\fontfamily{qpl}\selectfont 48 igi-gub \textit{\v{s}\`{a}} pa$_{5}$-sig}\\
   ``0;48 is the constant of a small canal''\\
   This suggest that the ratio of the depth of canal to that of its reserved water is $0;48=\frac{4}{5}$. In other words,    we   have 
\[ \frac{z'}{z}= \frac{4}{5}. \]  
Thus,   $z'=\frac{4}{5}z $ and the volume  of the reserved water should be 
$$V'=\frac{4}{5}V.$$
 So  it follows from \cref{equ-SMT24-b}   that
\begin{equation}\label{equ-SMT24-ba}
	V'= \frac{2}{5}xz(u+v).
\end{equation}

\subsubsection*{First Problem}\label{SSS-P1-SMT24} 
 Due to damage to the tablet, we are unable to establish with certainty the meanings of the technical expressions and calculations found in lines 1-25. We enumerate some of these ambiguities below: 

\vspace{2mm}
\noindent
\underline{Line 18: ``0;2 is the ratio $\cdots $''} 

\vspace{1mm}
\noindent
If we denote the width by $y$, this probably means $\frac{y}{x}=0;2$. In this case, according to lines 22-23, the value of $y$ is obtained by
\[ y=(0;2)x=(0;2)\times 45=1;30 \]
assuming that $x=45$. But this only adds to the confusion because similar terminologies (such as the upper breadth  and the lower breadth)    also occur in the text.  

\vspace{2mm}
\noindent
\underline{Line 22: ``15, the one which is to be added to the length''}

\vspace{1mm}
\noindent
If we assume that  $x=45$, this probably concerns the calculation
\[ x+15=45+15=1,0 \]
whose result is mentioned in line 25: ``Since 1,0 as the length  is said to you''.

\vspace{2mm}
\noindent
\underline{Lines 18-19:``$(0;30)\times (9;22,30)=4;41,45 $''}

\vspace{1mm}
\noindent
The result of this multiplication is called ``the area'' in lines 20-21. 

\vspace{2mm}
\noindent
\underline{Lines 20-21: ``$(1;30)\times (9;22,30)=14;3,45 $''}

\vspace{1mm}
\noindent
The result of this multiplication is called ``the false  area'' in lines 20-21. 

\vspace{2mm}
\noindent
\underline{Lines 22-24: ``$15\times (0;2)\times 3-0;30=(0;30)\times3 -0;30=1;30-0;30=1$''}

\vspace{1mm}
\noindent
The number 0;30 is called ``the one which is subtracted from the width'' in line 23.

\vspace{2mm}
\noindent
\underline{Line 24: ``$1\times (9;22,30)=9;22,30$''}

\vspace{1mm}
\noindent
We have been unable to reach a conclusion concerning this multiplication.

\vspace{2mm}
\noindent
\underline{Line 25:  ``$1,0+9;22,30=1,9;22,30$''}

\vspace{1mm}
\noindent
We have also ben unable to reach a conclusion concerning   this addition.

\begin{figure}[H]
	\centering
	\includegraphics[scale=1]{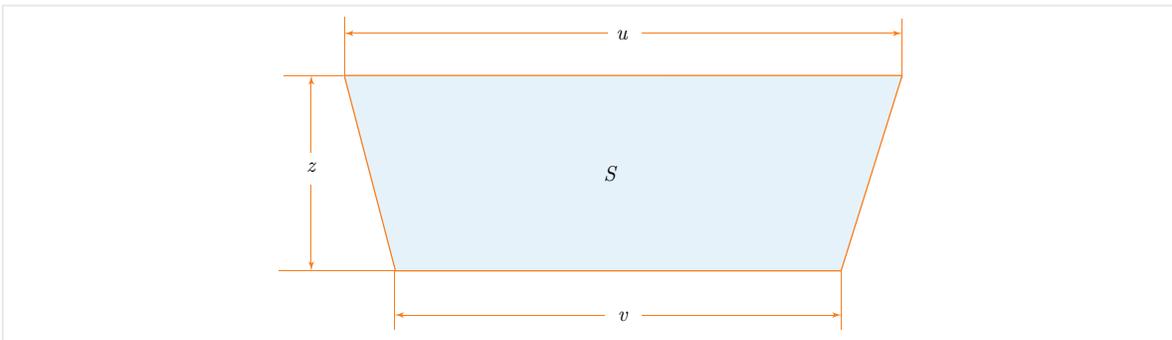}
	\caption{Cross-section of a trapezoidal canal}
	\label{Figure2}
\end{figure}
 
In spite of these uncertainties,   the scribe of this  text  is able to compute both   the dimensions of a canal  in general  and also solve  a quadratic equation  to find the upper breadth  of the canal.  We have shown the trapezoidal  cross-section  of a general trapezoidal  canal  in \cref{Figure2}.

Note that the calculations in lines 34-40 show  that the scribe is assuming the following relations between the dimensions of the canal:
\begin{equation}\label{equ-SMT24-c}
	\begin{dcases}
		v=\frac{u}{2}+0;30\\
		z=12\left(0;30+\frac{1}{12}(u-v)\right).
	\end{dcases}
\end{equation}

By looking carefully at the scribe's calculations in lines 26-33, it can be seen that he has solved the following quadratic equation:
\begin{equation}\label{equ-SMT24-d}
	(14;3,45) u^2-(1,9;22,30)u=4;41,15.
\end{equation}
This equation \cref{equ-SMT24-d} has been solved by the usual method of completing the squares-a standard method used by   Babylonian and Elamite scribes  to solve quadratic equations. This method was called {\fontfamily{qpl}\selectfont \textit{Tak\={\i}ltum}} in Babylonian texts  (see \cite{Mur03-1}, for a discussion on this topic).  

Now, let us use this method and solve the quadratic equation \cref{equ-SMT24-d} as follows:

\begin{align*}
	&~~  (14;3,45) u^2-(1,9;22,30)u=4;41,15 \\
	\Longrightarrow~~&~~  (14;3,45)\times (14;3,45) u^2-(14;3,45)\times(1,9;22,30)u=(14;3,45)\times(4;41,15) \\
	\Longrightarrow~~&~~   (14;3,45)^2 u^2-(1,9;22,30)\times \Big((14;3,45)u\Big)=1,5;55,4,41,15 \\  
	\Longrightarrow~~&~~  \Big((14;3,45)u\Big)^2-2\times (34;41,15)\times \Big((14;3,45)u\Big)=1,5;55,4,41,15 \\ 
	\Longrightarrow~~&~~   \Big((14;3,45)u\Big)^2-2\times (34;41,15)\times \Big((14;3,45)u\Big)+(34;41,15)^2\\ ~~&~~\hspace{8cm}  =1,5;55,4,41,15+(34;41,15)^2 \\ 
	\Longrightarrow~~&~~   \Big((14;3,45)u - 34;41,15 \Big)^2 =1,5;55,4,41,15+20,3;13,21,33,45 \\
	\Longrightarrow~~&~~   \Big((14;3,45)u - 34;41,15 \Big)^2 =21,9;8,26,15 \\ 
	\Longrightarrow~~&~~    (14;3,45)u - 34;41,15   =\sqrt{21,9;8,26,15} \\ 
	\Longrightarrow~~&~~    (14;3,45)u - 34;41,15   =\sqrt{(35;37,30)^2} \\ 
	\Longrightarrow~~&~~    (14;3,45)u - 34;41,15   =35;37,30 \\
	\Longrightarrow~~&~~    (14;3,45)u  =35;37,30+ 34;41,15   \\ 
	\Longrightarrow~~&~~    (14;3,45)u  =1,10;18,45  \\ 
	\Longrightarrow~~&~~    u  = \frac{1}{(14;3,45)} \times (1,10;18,45)\\ 
	\Longrightarrow~~&~~    u  = \left(0;4,16\right)\times (1,10;18,45)  \\ 
		\Longrightarrow~~&~~    u  = 5.              
\end{align*}
So
\begin{equation}\label{equ-SMT24-e}
	u=5.
\end{equation}
Now,  according to lines 33-36,  we can   find the values of $v$ and $z$. On the one hand, by \cref{equ-SMT24-c} and  \cref{equ-SMT24-e}, we have 
\begin{align*}
	v&=\frac{u}{2}+0;30\\ 
	&=\frac{5}{2}+0;30\\
	&=2;30+0;30\\
	&=3.
\end{align*}
Hence,
\begin{equation}\label{equ-SMT24-f}
	v=3.
\end{equation}
On the other hand, it follows from \cref{equ-SMT24-c} and  \cref{equ-SMT24-f}  that

\begin{align*}
	z&=12\left(0;30+\frac{1}{12}(u-v)\right)\\ 
	&=12\left(0;30+\frac{1}{12}(5-3)\right)\\
	&=12\left(0;30+\frac{2}{12}\right)\\
	&=12\times (0;30+0;10)\\
	&=12\times (0;40)
\end{align*}
which implies that
\begin{equation}\label{equ-SMT24-g}
	z=8.
\end{equation}

Finally, a verification process seems to begin at line 37. First,  the area\index{area of a trapezoid} of the trapezoidal\index{trapezoidal cross-section} cross-section\index{cross-section of a canal} $S$ is obtained by using \cref{equ-SMT24-a}, \cref{equ-SMT24-e}, \cref{equ-SMT24-f} and  \cref{equ-SMT24-g} as follows:
\begin{equation}\label{equ-SMT24-h}
	S =\frac{z(u+v)}{2} =\frac{8(5+3)}{2} =4\times 8=32.
\end{equation}
Then, the length\index{length} $x$ of the canal\index{canal} is obtained   by using the known volume\index{volume of a trapezoidal canal}  $V=24,0$ (mentioned in line 3). In fact, it follows from \cref{equ-SMT24-b} and  \cref{equ-SMT24-h} that
\begin{align*}
	x&=\frac{V}{S}\\ 
	&=\frac{24,0}{32}\\
	&= \frac{1}{32} \times (24,0)\\
	&=\left(0;1,52,30\right)\times (24,0)\\
	&=45.
\end{align*}

\subsubsection*{Second Problem}\label{SSS-P2-SMT24}
In the lines regarding the second problem, we can recognize several typical expressions found in excavation problems\index{excavation problems} as follows:

\vspace{2mm}
\noindent
(1) Reverse, line 2: ``I excavated''.\\
(2) Reverse, line 4: ``the area\index{area} of the large hole''.\\
(3) Reverse, line 5: ``the depth\index{depth (of a canal)} of the former canal\index{canal}''.\\
(4) Reverse, line 6: ``I cut off $\cdots $ for the small''.\\
(5) Reverse, line 8: ``the amount by which (the length\index{length} of) the (large) hole exceeded (the length\index{length} of) the (small) hole''.\\
(6) Reverse, line 10: ``the reciprocal of a number\index{reciprocal of a number} of 12 of the depth\index{depth (of a canal)}''.\\

\begin{figure}[H]
	\centering
	\includegraphics[scale=1]{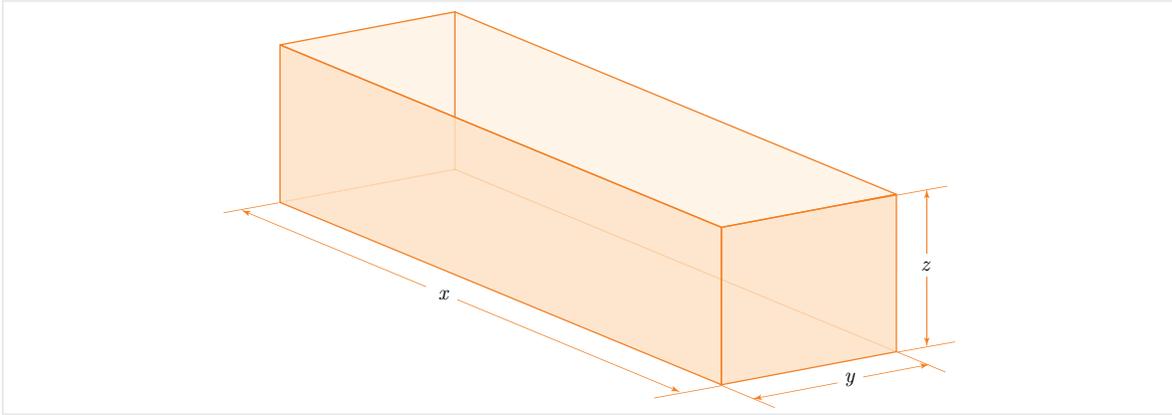}
	\caption{A canal with rectangular cross-section}
	\label{Figure3}
\end{figure}

Although we cannot restore the statement of this problem either, it seems that a canal\index{canal} has been enlarged, that is, the bottom of a canal\index{canal} has been deepened. Let $x$, $y$ and $z$ denote the length\index{length}, the width\index{width} and the depth\index{depth (of a canal)} of a canal\index{canal} whose cross-section\index{cross-section of a canal} is rectangular\index{rectangular cross-section} (see \cref{Figure3}).

Judging from the calculations performed in the text, the equations dealt with in this problem are:
\begin{equation}\label{equ-SMT24-i}
	\begin{dcases}
		x-y=0;10\\
		z=12(x-y)\\
		z(x^2+y^2) +xy(z+1)+\frac{1}{13}\Big(x^2+y^2\Big)=1;15.
	\end{dcases}
\end{equation}
Before discussing the geometrical significance of these equations and the nature of the hole ({\fontfamily{qpl}\selectfont \textit{kalakkum}}), we   analyze   the solution given in lines 7-24. 

According to line 7, we can multiply both sides of the third equation in \cref{equ-SMT24-i} by 13. Since $13\times (1;15)=16;15$, we  get
\begin{equation}\label{equ-SMT24-j}
	13z(x^2+y^2) +13xy(z+1)+ x^2+y^2 =16;15.
\end{equation}
At this point, the scribe appears to have used the algebraic identity\index{algebraic identity}
\begin{equation}\label{equ-SMT24-k}
	x^2+y^2=(x-y)^2+2xy.
\end{equation}
According to lines 8-9, since $x-y=0;10$, we can use \cref{equ-SMT24-j} and \cref{equ-SMT24-k} to write
\begin{align*}
	&~~  13z(x^2+y^2) +13xy(z+1)+ x^2+y^2 =16;15 \\
	\Longrightarrow~~&~~  13z\Big((x-y)^2+2xy\Big) +13xy(z+1)+ (x-y)^2+2xy =16;15 \\
	\Longrightarrow~~&~~   13z (x-y)^2+26xyz +13xyz+13xy+ (x-y)^2+2xy =16;15 \\ 
	\Longrightarrow~~&~~   13z (x-y)^2+(3\times 13)xyz + (13+ 2)xy =16;15-(x-y)^2  \\  
	\Longrightarrow~~&~~   13z (x-y)^2+(3\times 13)xyz + (13+ 2)xy =16;15-(0;10)^2  \\ 
	\Longrightarrow~~&~~   13z (x-y)^2+(3\times 13)xyz + (13+ 2)xy =16;15-0;1,40                       
\end{align*}
thus
\begin{equation}\label{equ-SMT24-l}
	13z (x-y)^2+(3\times 13)xyz + (13+ 2)xy =16;13,20.
\end{equation}
(Following the scribe's calculations, we did not substitute $x-y=0;10 $ in term  $13z (x-y)^2 $ to simplify \cref{equ-SMT24-l}  more.) 

  In lines 10-11, the scribe calculates the reciprocal of  \index{reciprocal of a number}   $z$. It follows from the first two equations in \cref{equ-SMT24-i} that
\begin{align*}
	&~~ z =12(x-y) \\
	\Longrightarrow~~&~~  \frac{1}{z}=\frac{1}{12}  \times  \frac{1}{x-y} \\
	\Longrightarrow~~&~~  \frac{1}{z}=\left(0;5\right) \times  \frac{1}{(0;10)} 
\end{align*}
which gives us
\begin{equation}\label{equ-SMT24-m}
	\frac{1}{z}=0;30.
\end{equation}
Next, according to lines 11-20, we multiply both sides of \cref{equ-SMT24-l} by $1/z$  and then use  \cref{equ-SMT24-l} to find the value of $xy$ as follows:
\begin{align*}
	&~~13z (x-y)^2+(3\times 13)xyz + (13+ 2)xy =16;13,20  \\
	\Longrightarrow~~&~~    \frac{1}{z} \times \Big(13z (x-y)^2+(3\times 13)xyz + (13+ 2)xy\Big) = \frac{1}{z} \times (16;13,20) \\
	\Longrightarrow~~&~~   13  (x-y)^2+(3\times 13)xy  + (13+ 2) \Big( \frac{1}{z} \times xy\Big)  = \frac{1}{z} \times (16;13,20)\\
	\Longrightarrow~~&~~   13  (x-y)^2+(3\times 13)xy  + \Big((13+ 2)\times \left(0;30\right)\Big)  xy  =\left(0;30\right)\times (16;13,20)\\
	\Longrightarrow~~&~~   13 \times  (0;10)^2+\Big((3\times 13)  +  13\times  \left(0;30\right)+2\times \left(0;30\right)\Big)  xy  =8;6,40\\
	\Longrightarrow~~&~~     13 \times  (0;1,40) + (39 +  6;30+1)   xy  =8;6,40\\
	\Longrightarrow~~&~~   0;21,40 +(46;30)  xy  =8;6,40\\
	\Longrightarrow~~&~~   (46;30)  xy  =8;6,40-0;21,40  \\
	\Longrightarrow~~&~~   (46;30) xy  =7;45   \\
	\Longrightarrow~~&~~     xy  = \frac{1}{(46;30)} \times (7;45)\\
	\Longrightarrow~~&~~     xy  = \frac{1}{6\times(7;45)} \times (7;45) 
\end{align*}
which implies that
\begin{equation}\label{equ-SMT24-n}
	xy=0;10.
\end{equation}

The last part of the text proceeds with the common Babylonian method (completing the square). According to lines 20-21, by \cref{equ-SMT24-i} and  \cref{equ-SMT24-n}, we can write: 

\begin{align*}
	\dfrac{x+y}{2}&=\sqrt{\left(\dfrac{x-y}{2}\right)^2+xy}\\
	&=\sqrt{\left(\dfrac{0;10}{2}\right)^2+0;10}\\
	&=\sqrt{\left(0;5\right)^2+0;10}\\
	&=\sqrt{0;0,25+0;10}\\
	&=\sqrt{0;10,25}\\
	&=\sqrt{(0;25)^2} 
\end{align*}
So
\begin{equation}\label{equ-SMT24-o}
	\dfrac{x+y}{2}=0;25.
\end{equation}
Now, we are in the position to compute the values of $x$ and $y$. According to lines 23-24, it follows from \cref{equ-SMT24-i} and  \cref{equ-SMT24-o} that
\begin{align*}
	x&=\dfrac{x+y}{2}+\dfrac{x-y}{2}\\
	&= 0;25+0;5\\
	&=0;30 
\end{align*}
and
\begin{align*}
	y&=\dfrac{x+y}{2}-\dfrac{x-y}{2}\\
	&= 0;25-0;5\\
	&=0;20.
\end{align*}
Therefore, the required values of the length\index{length} and the width\index{width}  are
\[ x=0;30\ \ \ \text{and}\ \ \ y=0;20. \]

\begin{figure}[H]
	\centering
	\includegraphics[scale=1]{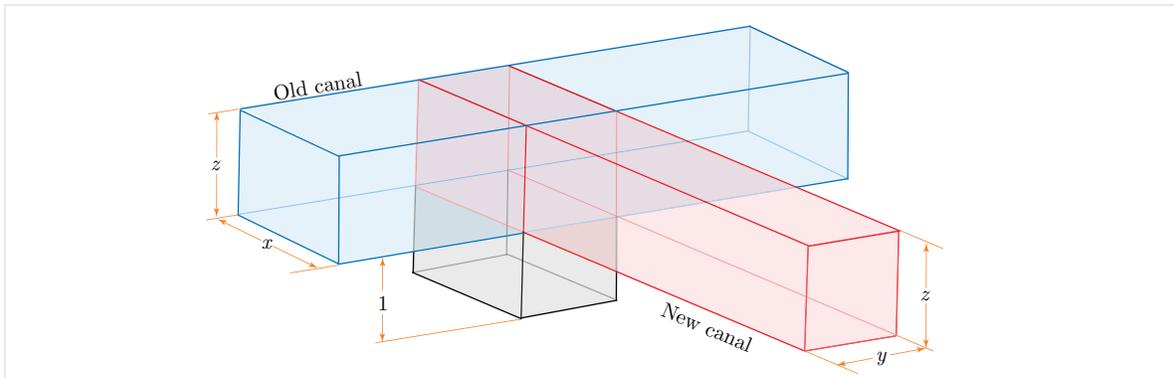}
	\caption{Intersection of  two canals}
	\label{Figure4}
\end{figure}

Although we can now understand the mathematical meaning of the second problem,   difficult questions  still remain: \\
What are really $x$ and $y$ in this problem?; and\\
 What is the relation between ``a canal\index{canal}'' ({\fontfamily{qpl}\selectfont \textit{atappum}}) and ``a hole'' ({\fontfamily{qpl}\selectfont \textit{kalakkum}})? \\
  To answer these questions, we may need to  consider a  hypothetical  situation   where two  canals intersect as shown in \cref{Figure4}. In this situation,      an old canal\index{canal} of  depth\index{depth (of a canal)} $z$ and   width\index{width} $x$ is joined by    a new canal\index{canal} of  depth\index{depth (of a canal)} $z$ and width\index{width} $y$ at right angles. The intersection between the old and new   canal\index{canal}s   is a rectangle\index{rectangle} with sides $x$ and $y$  and  deepened by 1 {\fontfamily{qpl}\selectfont k\`{u}\v{s}}\index{kuz@k\`{u}\v{s} (length unit)}, probably to  allow the deposit of silt\footnote{In modern Japan, people sometimes see the same device for  irrigation canal\index{irrigation canal}s\index{canal}.}.

\begin{figure}[H]
	\centering
	\includegraphics[scale=1]{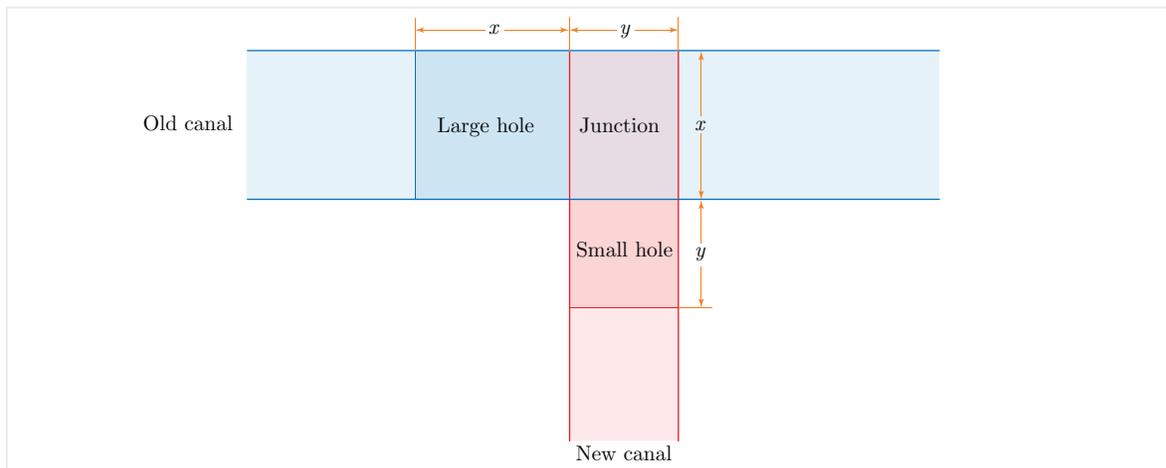}
	\caption{Top view of two intersecting canals and their junction}
	\label{Figure5}
\end{figure}

The top view of this intersection  is pictured  with related dimensions in  \cref{Figure5}. The rectangle in the east-west direction is the old canal and the one in the north-south direction is the new canal.  It is clear from the figure that the junction\index{intersection of two canals} of the two canal\index{canal}s is a cube of length\index{length} $x$, width\index{width} $y$ and depth\index{depth (of a canal)} $1+z$. According to the text, there are two  holes adjoining to this intersection\index{intersection of two canals}   which we have named the large hole and the small hole. Note that  the large hole ({\fontfamily{qpl}\selectfont \textit{kalakki} gal})   has dimensions $x$, $x$, $z$ and the small hole ({\fontfamily{qpl}\selectfont \textit{kalakki} tur}) has dimensions $y$, $y$, $z$. The unknown quantities asked for in this problem are the length\index{length} $x$  of the large hole     and the length\index{length} $y$ of the small hole.  One might expect here that $x$  is called, for example, {\fontfamily{qpl}\selectfont u\v{s} \textit{kalakki} gal} ``the length\index{length} of the large hole'' instead of  {\fontfamily{qpl}\selectfont \textit{kalakki} gal}  ``(the length) of the large hole''. Since $x$ is originally the width\index{width} of the old canal\index{canal}, the Susa scribe\index{Susa scribes} might have omitted {\fontfamily{qpl}\selectfont u\v{s}} ``the length\index{length}'' in order to avoid confusion. In fact, neither {\fontfamily{qpl}\selectfont u\v{s}} ``the length\index{length}''  nor {\fontfamily{qpl}\selectfont sag} ``the width\index{width}'' occurs in the text at all. 

Now, let us return to the third equation in \cref{equ-SMT24-i}. We   interpret the first and the second terms of the left-hand side, i.e., $z(x^2+y^2)=zx^2+zy^2 $ and $xy(z+1) $, as the sum  of the volumes\index{volume of a hole} of the large hole, the small hole and the junction\index{intersection of two canals} of the two canal\index{canal}s respectively. However, as to the third term $\frac{1}{13}(x^2+y^2) $, we are of the view  that it does not have any geometrical meaning and is   added to these volumes\index{volume of a hole}  to make the equation more complicated. There are examples of similar equations   in the so-called \textit{series texts}\index{series texts}, that is, the ancient drill books in mathematics. For example, the following system of equations  is given in \index{VAT 7537@\textbf{VAT 7537}} the text of  \textbf{VAT 7537}\footnote{Tablet \textbf{VAT 7537} is a Babylonian mathematical text in the Berlin Museum which originally was published by Neugebauer\index{Otto Neugebauer} in \cite{Neu35}. For more  information about this tablet and its text, see \cite{Mur91-1}.}:  
\begin{equation*}
	\begin{dcases}
		xy=10,0\\
		x^2+\frac{1}{11}\Bigg(2\left(\frac{1}{7}\Big((2y+3(x-y)\Big)\right)^2 +x^2\Bigg)=16,40.
	\end{dcases}
\end{equation*}

\subsection{SMT No.\,25}
As mentioned earlier, the problem in   \textbf{SMT No.\,25}   concerns the dimensions of a canal.

\subsubsection*{Transliteration}\label{SS-TI-SMT25}  

\begin{note1} 
	\underline{Reverse,  Lines 25-34}\\
	(L25)\hspace{0mm}  [$ \cdots $ $ \cdots $ $ \cdots $ \textit{me}]-\textit{e \v{s}\`{a}} gi\v{s}-g[i $ \cdots $] 2 \textit{\v{s}u\v{s}\v{s}\={a}r} (\v{S}\'{A}RxDI\v{S}) \textit{\v{s}\`{a}} pa$_5$ GIG(?)\\
	(L26)\hspace{0mm}  [$ \cdots $ $ \cdots $] $ \cdots $ [$ \cdots $] u\v{s} 5-ta-\`{a}m \textit{\v{s}\`{a}} 6 \textit{\v{s}}[\textit{\`{a}} \v{s}\'{a}r] $ \cdots $ -\textit{ma}
	\\
	(L27)\hspace{0mm}  i[gi-5] \textit{pu-\c{t}\'{u}-\'{u}r} 12 \textit{ta-mar} 12 \textit{a-na} 6 \textit{\v{s}\`{a}} \v{s}\'{a}r \textit{i-\v{s}\'{i}-ma}\\
	(L28)\hspace{0mm}  1,12 \textit{ta-mar} 1 \v{s}\'{a}r \textit{\`{u}} 12 \textit{\v{s}u-\v{s}i mu-\'{u} a-na} 40 er\'{i}n-hi-a [gar]
	\\\
	(L29)\hspace{0mm}  igi-40 erín \textit{pu-\c{t}\'{u}-\'{u}r} 1,30 \textit{ta-mar} 1,30 \textit{a-na} 1,12 \textit{i-\v{s}\'{i}-ma}
	\\
	(L30)\hspace{0mm}  1,48 \textit{ta-mar} \textit{mu-\'{u} \v{s}\`{a}} er\'{i}n 1-k[am u\v{s}] 1 nindan 30 dagal
	\\
	(L31)\hspace{0mm}  [1],48 \textit{mu-\'{u} \v{s}u-up-lu mi-nu} igi-48 igi-gub pa$_5$ \textit{pu-\c{t}\'{u}}-[\textit{\'{u}r}]
	\\
	(L32)\hspace{0mm}  [1,15 \textit{t}]\textit{a-mar} 1,15 \textit{a-na} 1,48 \textit{me-e \v{s}\`{a}} [er\'{i}n 1-kam \textit{i-\v{s}\'{i}-ma}]\\
	(L33)\hspace{0mm} [2,15 \textit{t}]\textit{a-mar} igi-30 dagal \textit{pu-\c{t}\'{u}-\'{u}r} 2 \textit{ta-mar}
	\\
	(L34)\hspace{0mm}  [2,15 \textit{a-n}]\textit{a} 2 \textit{i-\v{s}\'{i}-ma} 4,30 \textit{ta-mar} 4,30 \textit{\v{s}u-up-lu}
\end{note1}

\subsubsection*{Translation}\label{SS-TR-SMT25}  

\underline{Reverse,  Lines 25-34}
\begin{tabbing}
	\hspace{15mm} \= \kill 
	(L25)\> \tabfill{$ \cdots $ $ \cdots $ water of reed bed  $ \cdots $. 2,0 saros (of water) of a canal. $ \cdots $.}\index{reed bed}\index{canal} \\
	(L26)\> \tabfill{$ \cdots $ $ \cdots $ the length 5 ({\fontfamily{qpl}\selectfont nindan}\index{nindan (length unit)}) each that of 6 saros $ \cdots $ $ \cdots $.}\index{length}\\
	(L27)\> \tabfill{Make the reciprocal of 5, (and) you see 0;12. Multiply 0;12 by 6 saros, and}\index{reciprocal of a number} \\
	(L28)\> \tabfill{you see 1,12,0. 1 saros and 12 sixties is the (volume of) water. Put (it) down for 40,0 workers.}\index{volume of a canal}\\
	(L29)\> \tabfill{Make the reciprocal of 40,0 of the workers, (and) you see 0;0,1,30. Multiply 0;0,1,30 by 1,12,0, and}\index{reciprocal of a number} \\
	(L30)\> \tabfill{you see 1;48. (This is the volume of) the water (per 1 {\fontfamily{qpl}\selectfont nindan}\index{nindan (length unit)} in length) of the first(sic) worker. (If) the length is 1 {\fontfamily{qpl}\selectfont nindan}\index{nindan (length unit)}, the breadth is 0;30 ({\fontfamily{qpl}\selectfont nindan}\index{nindan (length unit)}), (and)}\index{volume of a canal}\index{breadth (of a canal)}\index{length}\\
	(L31)\> \tabfill{the (volume of) water is 1;48, what is the depth? Make the reciprocal of 0;48, the constant of a canal, (and)}\index{canal}\index{reciprocal of a number}\index{depth (of a canal)}\index{volume of a canal}\\
	(L32)\> \tabfill{you see 1;15. Multiply 1;15 by 1;48 of the water of a worker, and}\\
	(L33)\> \tabfill{you see 2;15. Make the reciprocal of 0;30 of the breadth, (and) you see 2.}\index{reciprocal of a number}\index{breadth (of a canal)} \\
	(L34)\> \tabfill{Multiply 2;15 by 2, and you see 4;30. 4;30 ({\fontfamily{qpl}\selectfont k\`{u}\v{s}}) is the depth.}\index{kuz@k\`{u}\v{s} (length unit)}\index{depth (of a canal)} 
	
\end{tabbing}

\subsubsection*{Mathematical Calculations}\label{SS-MC-SMT25}
The statement of the problem is almost entirely broken, and we have been unable to locate a similar problem to assist us in restoring it.  However, the text calculates the depth\index{depth (of a canal)} of a canal\index{canal} using   several conditions which may well have been provided in the statement. Judging from the calculations performed in the text, the missing conditions might have  been as follows:

\vspace{2pt}
\noindent
\underline{\textbf{Problem:}} Two thousand   four hundred workers (40,0 {\fontfamily{qpl}\selectfont er\'{i}n-hi-a}) built a canal\index{canal} whose width\index{width} is 0;30 {\fontfamily{qpl}\selectfont nindan}\index{nindan (length unit)}, and whose reserved water\index{depth (of the reserved water)} is 6 {\fontfamily{qpl}\selectfont \v{s}\'{a}r}\footnote{Note that  1 {\fontfamily{qpl}\selectfont \v{s}\'{a}r}=1 saros \index{saros} equals to 1,0,0 in sexagesimal\index{sexagesimal number system} number system.} (that is, 6,0,0 {\fontfamily{qpl}\selectfont volume-sar}\index{volume-sar (capacity unit)}). The workers were each assigned to dig a part of the canal\index{canal}   5 {\fontfamily{qpl}\selectfont nindan}\index{nindan (length unit)} in length\index{length}. What is the depth\index{depth (of a canal)} of the canal\index{canal}?

\vspace{2pt}

\begin{figure}[H]
	\centering
	\includegraphics[scale=1]{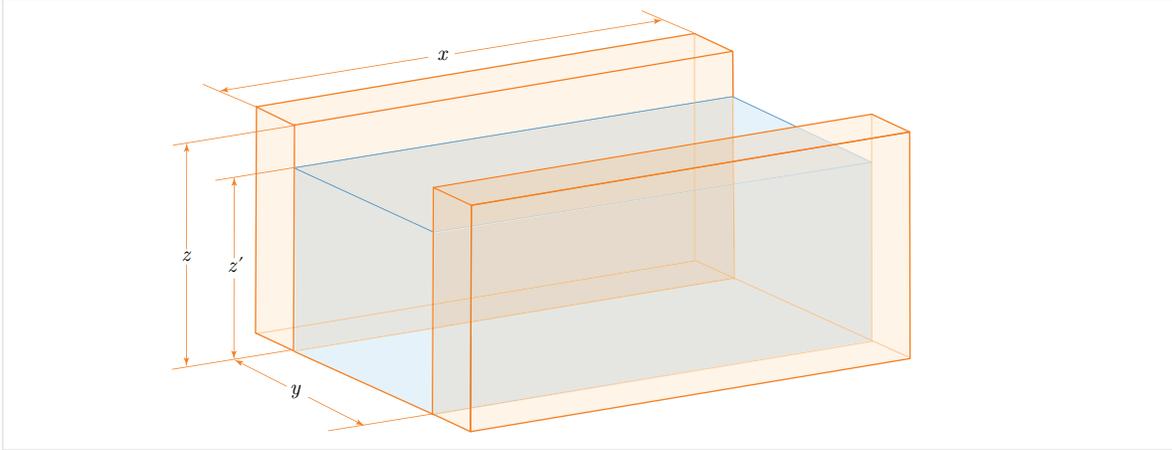}
	\caption{A canal and its reserved water}
	\label{Figure6}
\end{figure}

Consider a canal\index{canal} of  length\index{length} $x$,  width\index{width} $y$, and depth\index{depth (of a canal)} $z$ and denote  the depth\index{depth (of the reserved water)} of its reserved water\index{reserved water of a canal}    by  $z'$  as is shown in \cref{Figure6}. As we said before, we must have
\begin{equation}\label{equ-SMT25-a}
	\frac{z'}{z}=0;48.
\end{equation}

Let $V$ be the volume\index{volume of a canal} of the canal\index{canal}, $V'$ the volume\index{volume of a canal} of its  reserved water\index{reserved water of a canal},  $S$ the area\index{area of a rectangle} of its cross-section\index{cross-section of a canal} and $S'$ the area\index{area of a rectangle} of a part of the cross-section\index{cross-section of a canal} submerged in water. It is clear from the figure that
\begin{equation}\label{equ-SMT25-b}
	\begin{cases}
		V=xyz\\
		V'=xyz'\\
		S=yz\\
		S'=yz'.
	\end{cases}
\end{equation} 
Note that it follows from \cref{equ-SMT25-a} and  \cref{equ-SMT25-b}   that 
\begin{equation}\label{equ-SMT25-c}
	\frac{V'}{V}=\frac{S'}{S}=0;48.
\end{equation}

\begin{figure}[H]
	\centering
	\includegraphics[scale=1]{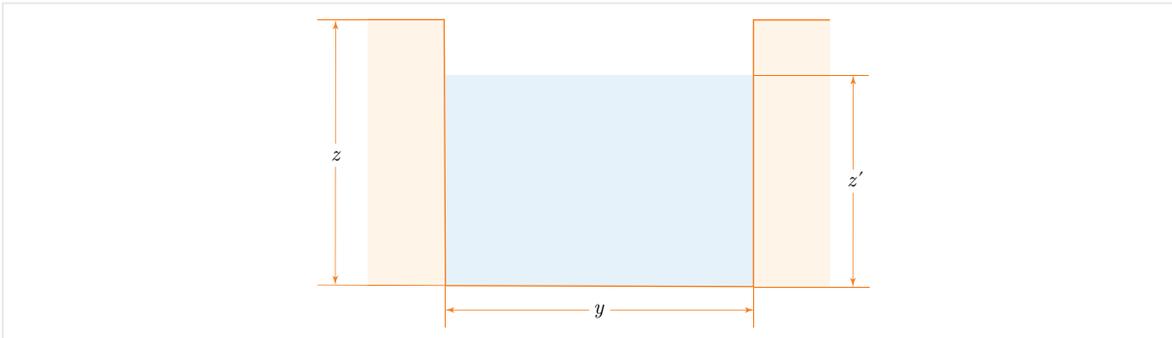}
	\caption{The cross-section of a canal and the level of its reserved water}
	\label{Figure7}
\end{figure}

In \cref{Figure7}, we have shown the cross-section  of the canal  and   its reserved water. The height of canal is $z'$ and that of reserved water is $z$.

In lines 26-30, the scribe has calculated  the volume\index{volume of a canal}   of the water per 1 {\fontfamily{qpl}\selectfont nindan}\index{nindan (length unit)} and per one worker, say $V'_1$. To do so, he first obtains the the volume\index{volume of a canal} $V'_0$  of water per 1 {\fontfamily{qpl}\selectfont nindan}\index{nindan (length unit)} and all workers as follows:
\begin{align*}
	V'_0&=\frac{V'}{5}\\
	&= \frac{1}{5} \times (6\ \text{{\fontfamily{qpl}\selectfont šár}}) \\
	&= (0;12)\times (6,0,0) \\
	&= 1,12,0 \\
	&= 1\ \text{{\fontfamily{qpl}\selectfont \v{s}\'{a}r}}\ 12\   \text{{\fontfamily{qpl}\selectfont \textit{\v{s}\={u}\v{s}i}}}. 
\end{align*}
(Note that 1 {\fontfamily{qpl}\selectfont \textit{\v{s}\={u}\v{s}i}} is equal to  $1,0=60$.)  
Then, he divides this volume\index{volume of a canal}  $V'_0=1,12,0 $ by the number of all workers to find $V'_1$, which   is also  called ``the water of the first worker'' in line 30:
\begin{align*}
	V'_1&=\frac{V'_0}{(40,0)}\\
	&= \frac{1}{(40,0)} \times (1,12,0) \\
	&= (0;0,1,30) \times (1,12,0) \\
	&= 1;48. 
\end{align*}
Note that  the value of $V'_1=1;48$ is the volume\index{volume of a canal} of the reserved water\index{reserved water of a canal} of the part of canal\index{canal} with length\index{length} 1.   That is,    $V'_1$ is obtained by assuming $x=1$ in the second equation of \cref{equ-SMT25-b}:  
\[ V'_1=1yz'=yz'=S'. \] 
This means $ V'_1$  is   equal to  the area\index{area of a rectangle} of the cross-section\index{cross-section of a canal} of the lower part of the canal\index{canal} submerged in  water: 
\begin{equation}\label{equ-SMT25-d}
	S'=1;48.
\end{equation} 
From \cref{equ-SMT25-b} and  \cref{equ-SMT25-d}, (according to lines 31-33) we can obtain $S$:
\begin{align*}
	S&=  \frac{S'}{(0;48)}\\
	&=\frac{1}{(0;48)}\times (1;48)\\
	&=(1;15)\times (1;48)
\end{align*}
which implies that
\begin{equation}\label{equ-SMT25-e}
	S=2;15.
\end{equation} 
Since  by the assumption $y=0;30$, we can (according to lines 33-34) use \cref{equ-SMT25-b} and  \cref{equ-SMT25-e} to find the depth\index{depth (of a canal)} of the canal\index{canal}, i.e., $z$,  as follows:
\begin{align*}
	z =  \frac{S}{y} 
	 =\frac{1}{(0;30)}\times (2;15) 
	 =2\times (2;15) 
	  =4;30.
\end{align*}
Therefore the depth\index{depth (of a canal)} of the canal\index{canal} is 
\begin{equation}\label{equ-SMT25-f}
	z=4;30.
\end{equation} 
Note that the depth\index{depth (of a canal)} of water can be easily computed by using  \cref{equ-SMT25-a} and  \cref{equ-SMT25-f} as follows:
\[ z'=(0;48)\times (4;30)=3;36. \]

\section{Applications of Mathematics  in Ancient Elam}
In this section, we consider the possible practical applications of mathematics    both in construction projects and in Elamite architecture. The reader will note  that   the mathematical skills demonstrated in  \textbf{SMT No.\,24}  and  \textbf{SMT No.\,25}   may well have been of considerable assistance to those charged with the construction of Elamite infrastructure and other substantial buildings at the behest of the Elamite rulers of the time.

\subsection{Canals}
\begin{figure}[H]
	\centering
	\includegraphics[scale=1]{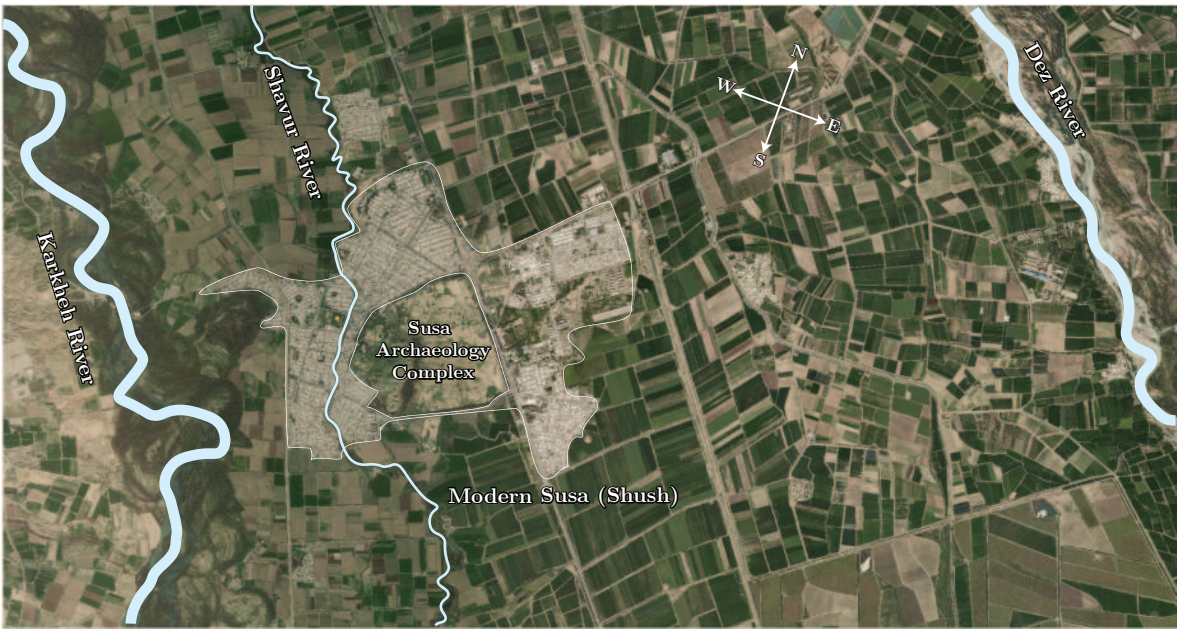}
	\caption{A satellite view of modern Susa (Google Map)}
	\label{SusaCity1}
\end{figure}

As for  any   ancient civilization,     water  played a decisive role in where the people of ancient Elam decided to settle. The Elamite people  were among the   very first farmers in the Ancient Near East   (see \cite{ABW18}), utilizing the nearby rivers    to irrigate their  lands and farms. For example, the ancient capital city of Susa was founded  in a region watered by at least two rivers\footnote{Professor Daniel Potts has made a strong case that there was only one river at Susa and that the modern Shavur river was the ancient course of the Karkheh river. See \cite{Pot99} for more details.}: the Karkheh and the Dez   (see \cref{SusaCity1}).

\begin{figure}[H]
	\centering
	\includegraphics[scale=1]{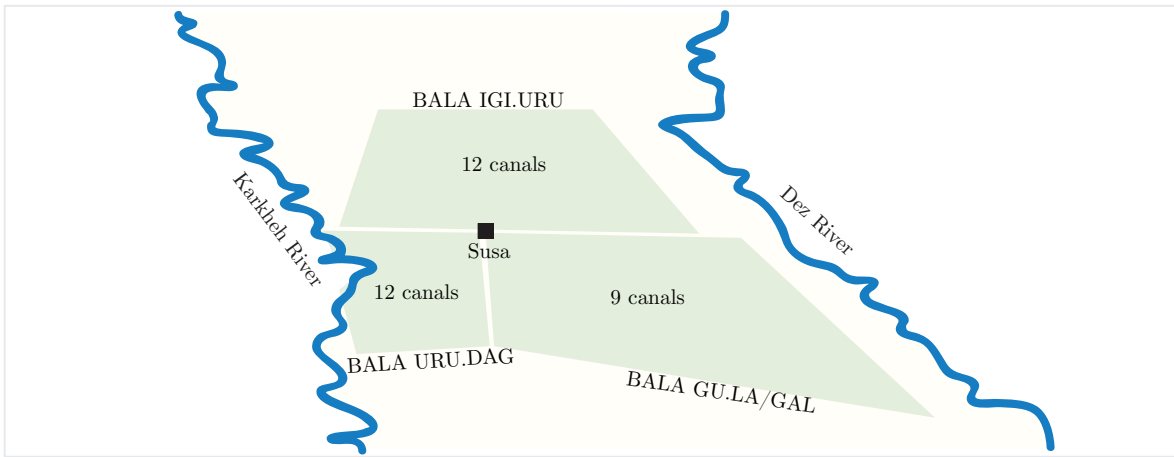}
	\caption{Possible distribution  of   irrigation canals near  Susa (adapted from \cite{Gra18})}
	\label{Susacanal-1}
\end{figure}

Although the course of the Karkheh in ancient Elam has been the subject of considerable scholarly debate, it has been suggested that since the Karkheh river is at a higher level than   the fields   between the Karkheh and the Dez, the people of Susa may  have   used the   Karkheh's water to irrigate the   land  around the city    or save it for   future purposes (probably for drinking or making mud bricks).    According to some research (see  \cite{Gra18,Sad14}), there were at least four agricultural sections  near   ancient Susa   which contained  more than  40 small irrigation canals. Based on the  analysis of   other legal, economic and administrative clay tablets excavated from Susa, scholars have been  able to determine the names and  locations of many of these  canals in each of these   agricultural sections. These tablets date to the Sukkalmah Dynasty 1900-1500 BC, as do the \textbf{SMT}. We have shown   three of these sections   and the  number of canals they contained   in \cref{Susacanal-1}. The location of the fourth section, described as being ``on the other bank of the \textit{Zamun}'', cannot be established as it is not known which river was at that time known as the \textit{Zamun}.

In addition to this network of irrigation canals around ancient Susa, there   was  also     the  ancient Harmushi canal system, which  consisted of  a group  of  connected   canals watering the vast area of  lands between the Karkheh and the Dez river.    Archaeological evidence  show that  some    branches  of this network  date  to the Middle Elamite period  (see \cite{Ada62, Ras20, Ras21}).

\begin{figure}[H]
	\centering
	\includegraphics[scale=1]{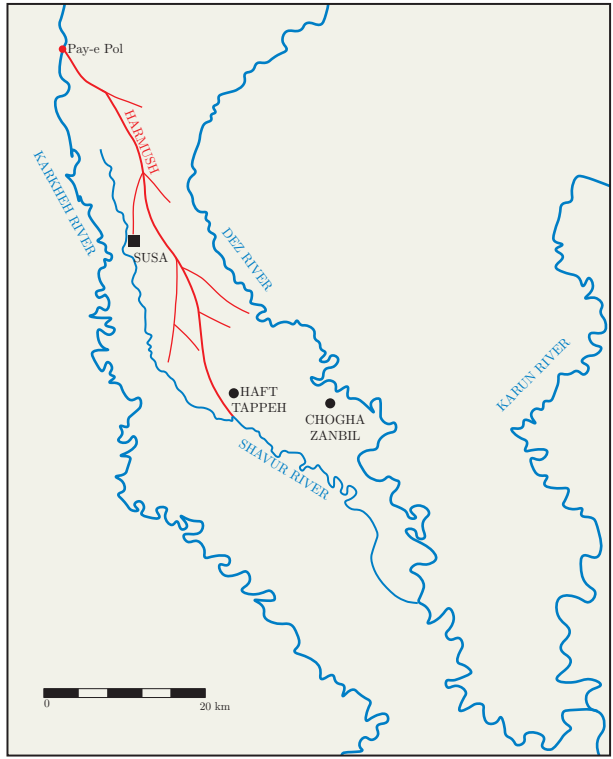}
	\caption{The hypothetical course of the Harmushi canal network}
	\label{HarmushiCourse}
\end{figure}

The main branch of this canal system was called the \textit{Harmushi} and connected to the Karkheh river at the point called \textit{Pay-e Pol} (literally, ``at the foot of the bridge''), where there  are still the ruins of an ancient dam which   regulated the river into different branches. The length of this canal is thought  to have been  approximately $ 50 km $ in which case  it could have  irrigated  nearly $200 km^2$ of the land    between the Karkheh and the Dez rivers. \cref{HarmushiCourse} shows a hypothetical course of the Harmushi canal network which has been adopted from Adams \cite[Fig. 4]{Ada62}.  The Harmushi system  had several subbranches which is thought to  provide  water  for   the livestock  of        nomadic tribes.    Many  western travelers and scholars, who  visited Khuzistan during last two centuries,    mentioned   the Harmushi canal in their reports  and asserted that it was the main source for  irrigating a large area    around Susa (see \cite{Ehl75,Fie39, Lor08,Raw36}). Although this ancient canal system   had been  irrigating  the  Susa area for more than 3000 years,  it was finally  replaced by the modern canal network after the construction of  a  concrete dam built on the Dez river in 1963.

 To understand the historical significance of this canal and its importance to the local population,   even  though there is almost no physical trace of the original system, the name of Harmushi still lives on in the surname of many   people whose  ancestors lived alongside it   for centuries.\footnote{The complete name of the first author used to be ``Nasser Heydari Harmushi''.}

\begin{figure}[H]
	\centering
	\includegraphics[scale=1]{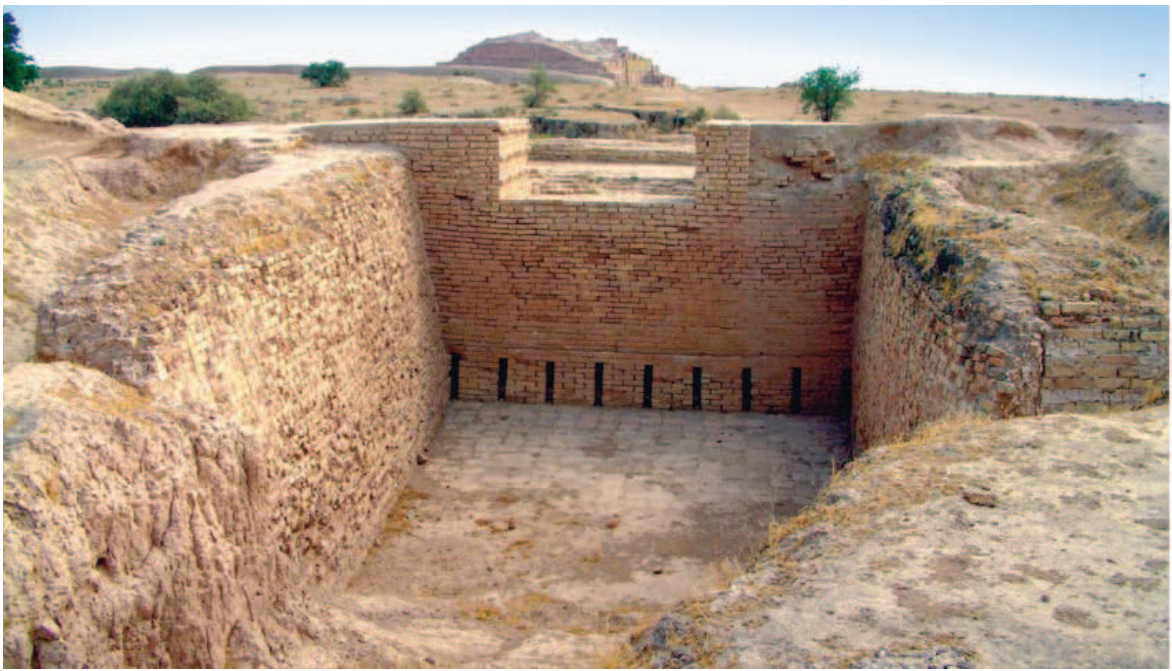}
	\caption{A water reservoir at the holy city of Chogha Zanbil (Credit: The World Heritage of Chogha Zanbil)}
	\label{Choghazanbil1}
\end{figure}

In addition to  canals, there are archaeological sites  near modern Susa  containing the remains of  other water structures one of which is located in   the   Elamite holy city of   Chogha Zanbil\footnote{The holy city of Chogha Zanbil (also known as \textit{Dur-Untash}, or City of Untash,  in Elam) located at $ 38km $ to the southeast of Susa  was founded by the Elamite king \textit{Untash-Napirisha}  (1275--1240 BCE) to be the religious center of Elam. The principal element of this complex is an enormous ziggurat  dedicated to the  Elamite divinities Inshushinak   and   Napirisha. For more information, see \url{https://whc.unesco.org/en/list/113/} or \url{http://tchoghazanbil.com/}.}. There is     a hydraulic complex    consisting of   canals  and a water reservoir  of length  $10.70m$, width  $7.25m$ and  depth  $4.35m$ whose capacity is about $332m^3$ (see  \cref{Choghazanbil1}).

Although as yet there is no satisfactory explanation for how the hydraulic mechanism of this sophisticated piece of engineering  worked, scholars have suggested different hypotheses. Some scholars including Ghirshman  believed that this reservoir  was fed by the Karkheh  river via   the Harmushi canal, which was specially built to supply  the new holy city (see \cite{Ghi68, Ste67}). Drinking water was  then thought to be distributed throughout the city by a system of smaller canals.

While Ghirshman's opinion prevailed for decades,  modern archaeological  research has raised serious doubts about his analysis (see \cite{Saf19}).    Mofidi-Nasrabadi in \cite{Mof07} has suggested  that this water reservoir  is actually a part of a drainage system   devised by the Elamite engineers to remove   floodwater from the holy city during the rainy season.

In either case,   the mathematical knowledge of Elamite scribes  might have been  of assistance in the construction of   such  a complex   structure. Their mathematical skills could have been applied to estimate the amount of time and number of workers needed for such a large engineering project.

\section{Conclusion}
The   interpretation of  problems contained  in   \textbf{SMT No.\,24} and \textbf{SMT No.\,25} taken together with the archaeological research confirming the existence of a substantial canal network at the time the problems were inscribed on the tablet, strongly infer that the Susa scribes were interested in the practical application   of mathematics to the challenges faced by those living alongside them.  These texts suggest that not only was     mathematics   taught   in a ``Susa School of Mathematics'', but also that the   scribes  used  their mathematical skills to address issues arising in the design and construction of   structures, such as canals, which facilitated both the   agricultural and economic  development of ancient Elam.

\end{document}